\documentclass{article} 
\usepackage[utf8]{inputenc} \usepackage{amsfonts} \usepackage{amsmath} \usepackage{indentfirst}
\usepackage{biblatex}
\newcommand{\isint}{\in \mathbb{Z}^{+}}
\newcommand{\zzp}{\mathbb{Z}^{+}}

\newcommand{\floor}[1]{\left \lfloor #1 \right \rfloor}
\newcounter{count}
\newcounter{cases}
\newcommand{\boldstart}[2]{\noindent\textbf{#1: }#2}

\newcommand{\lemma}[2]{\boldstart{Lemma \stepcounter{count}\thecount}{#1}\\\boldstart{Proof of Lemma \thecount}{#2}\\} 
\newcommand{\corollary}[2]{\boldstart{Corollary \stepcounter{count}\thecount}{#1}\\\boldstart{Proof of Corollary \thecount}{#2}\\} 
\newcommand{\induction}[2]{\boldstart{Base case(s)}{#1}\\\boldstart{Inductive Step}{#2}\\\\} 
\newcommand{\pair}[2]{\left< #1 , #2 \right>}

\title{Using Elementary Techniques to Characterize the Relationship Between Wythoff's Game and the Golden Ratio}
\author{Vincent Wang \and Nikhil Sampath \and Eric Yule\thanks{Special thanks} \and Ethan Wang\footnotemark[1]}
\date{January 2023}

\begin{document}

\maketitle

\section{Introduction}
Wythoff's game is a modification of the well-known game of ``nim." Wythoff's game, which does not resemble the Fibonacci sequence, has direct relation to the Golden ratio. We will explore the sequence behind this surprising relationship, and consider the implications of our elementary methods.

Wythoff's game is also known as the game of Ps and Qs. The game has two players taking turns removing chips of two types: P and Q (hence the game's name). Each turn has a player take any number of available Ps (and only Ps), any number of available Qs (and only Qs), or any same number of available Ps and Qs (for example a turn can have 2 Ps and 2 Qs being taken). Players alternate their turns until the player to remove the last chip wins \cite{wikipedia}.

Let us begin with notations and conventions. The first player is the player who took the first turn. The second player is the other player. Denote the state of $x$ Ps and $y$ Qs by $\pair{x}{y}$. If we start with $\pair{1}{2}$, no matter what turn the first player took, the second player can take a turn to win. We call a state where the second player can force a win a losing state, so that $\pair{1}{2}$ is losing. Analogously define a state where the first player can force a win a winning state.

We have some early observations. We observe that ``most" states are winning states, which motivates us to find the ``fewer" losing states. With some work, we find that $\pair{3}{5}$ is the next losing state (because of symmetry of the game with respect to P and Q, take $x < y$). Our final observation is that if $\pair{x}{y}$ has the property that exactly one of $x$ or $y$ was a component of a previously found losing state $\pair{z}{w}$, we can sometimes reduce $\pair{x}{y}$ to $\pair{z}{w}$, making $\pair{x}{y}$ a winning state. Wythoff, the mathematician for which the game is named after, characterizes this precisely in his fourth observation \cite{wythoff}. Our observations and computations\footnote{For reference Wythoff lists the first few pairs of losing states on page 200 of \cite{wythoff}} of  the first few losing states motivate the following conjecture:

Define the sequences $p(n)$, $q(n)$ for $n \geq 1$ recursively as follows: $p(1) = 1$, $q(n) = p(n) + n$ for all $n \geq 1$, and $p(n + 1)$ the smallest positive integer not contained in \[ U_n = \{ p(1), p(2), p(3), \dots , p(n), q(1), q(2), q(3), \dots , q(n) \}. \] Then $\pair{p(n)}{q(n)}$ for $n \geq 1$ characterizes all losing states.

We believe some scholarship already shows that $p$ and $q$ solve the game. The conjecture is provable by inducting on the value of $y - x$ over all such states $\pair{x}{y}$ in addition to significant casework on potential moves by each player. Such casework and induction is rather uninteresting. We are not focused on the game theory, but rather the analysis of the sequences $p$ and $q$.  Hence we will not prove the previous conjecture in this paper. Upon further inspection of our first few losing states, we find that the ratio of $q(n)$ to $p(n)$ does not grow as $n$ grows. If we plot $q$ versus $p$ in the Cartesian coordinate plane we find that the set of all points $(p, q)$ closely follows the line through the origin with slope $\phi$, the Golden ratio. Consequently, we instead focus on the following conjecture: \[ p(n) = \floor{\phi n} \text{ and } q(n) = \floor{\phi ^2 n} \text{ for all } n \geq 1.\] This direct relationship between the Golden ratio $\phi$ and our sequences $p$ and $q$ is very motivating.

Although Wythoff proves our conjecture brilliantly, we will offer an alternate partial proof. Wythoff's solution is effective as it quickly exploits the properties of $p$ and $q$, but it is also contrite as these properties are used in unison. Our elementary methods are incremental and offer ideas that can be applied to other problems.

\section{Understanding $p$ and $q$}
We will understand the behavior of $p$ and $q$ before we formalize the relationship between $p$, $q$, and $\phi$. For convenience we have the following naming conventions. $p$ by itself refers to the sequence $p(1), p(2), p(3), \dots$. If $N$ is a $p$, $N$ is contained in $p$. $N$ $p$s refers to a collection of $N$ elements of $p$. We have analogous names for $q$. Denote the set $\{1, 2, 3, \dots , m\}$ by $[m]$. We begin with a couple lemmas. \newline

\lemma{$p$ is strictly increasing}{We show the case of $p(n + 2) > p(n + 1)$ for all $n \geq 1$. The case of $p(2) > p(1)$ is left to the reader. Then it remains to show $[p(n + 1)] \subseteq U_{n + 1}$. Note for $n \geq 1$, $[p(n + 1) - 1] \subseteq U_n$ by the minimality of $p(n + 1)$. Therefore \[ [p(n + 1) - 1] \cup \{p(n + 1)\} \subseteq U_n \cup \{p(n + 1)\} \subseteq U_n \cup \{p(n + 1), q(n + 1)\},\] or $[p(n + 1)] \subseteq U_{n + 1}$, as desired.}

\corollary{There cannot be two consecutive $q$s. Define $\Delta q(n) = q(n + 1) - q(n)$ for $n \geq 1$. Then $\Delta q(n) \geq 2$.}{Follows from Lemma 1 and $q(n) = p(n) + n$.}

\lemma{$p$ and $q$ partition $\zzp$. That is, given any $N \geq 1$, $N$ appears exactly once in $p$ or $q$ (but not both).}{This fact should be rather natural. After some experimentation, one should see that $q$ is generated ``first", while $p$ fills in the gaps left by $q$. For completeness, our two part proof is as follows. We first show that any $N \geq 1$ appears in some $U$, then show the elements of $U$ are distinct. As seen in the proof for $p$ increasing (Lemma 1), we have $[p(n + 1) \subseteq U_{n + 1}]$. Take $n$ sufficiently large, so $N \leq p(n + 1)$ (we can do this because $p$ is strictly increasing). Then $N \in [p(n + 1)]$, and $N \in U_{n + 1}$. $p$s are pairwise distinct as $p$ is increasing. $q$s are pairwise distinct as $q$ is increasing. It remains for us to show a $p$ cannot be a $q$. Consider any $p(k)$ and $q(j)$ for $k, j \geq 1$. If $j < k$, $p(k) \neq q(j)$ by construction. If $j \geq k$, $q(j) \geq q(k) = k + p(k) > p(k)$, completing the proof.}

\lemma{Similarly define $\Delta p(n) = p(n + 1) - p(n)$ for $n \geq 1$. Then $\Delta p(n) \in \{1, 2\}$.}{$p(n)$ is increasing, so then next $p$ is at least $p(n) + 1$. If $p(n) + 1$ is a $p$, we are done and $p(n + 1) = p(n) + 1$. Otherwise, $p(n) + 1$ is a $q$ as $p, q$ partition $\zzp$. There cannot be two consecutive $q$s by Corollary 2, so $p(n) + 2$ cannot be a $q$, and hence must be a $p$, again by $p, q$ being a partition of $\zzp$. In this case, $p(n + 1) = p(n) + 2.$ We have exhausted all cases, proving the lemma.}

\corollary{$\Delta q(n) \in \{2, 3\}$ for all $n \geq 1$}{Follows from $\Delta p(n) \in \{1, 2\}$ and $q(n) = p(n) + n$.}

\corollary{There do not exist three consecutive $p$s.}{Follows from the previous corollary, $q(1) = 2$, and $p$, $q$ partition $\zzp$.}

After much consideration, we realize that there is little more we can say directly about $p$ and $q$. We have reached a ``bottleneck." To push through this bottleneck, we have a clever construction, which is a key focus of our paper. \\

\lemma{$q(n) = p(p(n)) + 1$ for all $n \geq 1$.}{Consider the set $[q(n)]$. The $q$s of this set are precisely $\{q(1), q(2), q(3), \dots , q(n)\}$, and there are $n$ of them. $p$ and $q$ partition $\zzp$, so there are $q(n) - n = p(n)$ $p$s of $[q(n)]$ (1). There cannot be two consecutive $q$s, so $q(n) - 1$ is a $p$. Moreover, $q(n) - 1$ is the largest $p$ in $[q(n)]$ (2). Combining (1) and (2) gives $q(n) - 1 = p(p(n))$, as desired.}

The crux idea of this proof is that $p$ and $q$ partition $\zzp$, so each positive integer is either a $p$ or a $q$. This seemingly trivial idea of being one or the other is very powerful. More generally, the fact that for any item $x$ and set $S$, we have either $x \in S$ or $x \not \in S$. We believe that strong theorems can be proven using set constructions exploiting the properties of the elements of that set. Here we used the property that there cannot be two consecutive $q$s and $q(n) = p(n) + n$. Let us see the implications of this result. \\

\lemma{$\Delta p(n) = 2$ if and only if $n$ is a $p$ for all $n \geq 1$.}{We present a non-typical method of proof, which we verify first as follows. Let $a$, $b$, and $c$ be statements and $\neg a$ denote the negation of statement $a$. We claim that if we show $a \rightarrow b$, $a \rightarrow c$, $\neg a \rightarrow \neg b$, and $\neg a \rightarrow \neg c$, then $b \leftrightarrow c$. The reason for the truth of this is as follows. $\neg a \rightarrow \neg b$ and $\neg a \rightarrow \neg c$ become $b \rightarrow a$ and $c \rightarrow a$ as any statement is equivalent to its contrapositive. From this we find $a \leftrightarrow b$ and $a \leftrightarrow c$, so $b \leftrightarrow c$, as desired. 

In our proof, $a$ is ``$p(n) + 1$ is a $p$", $b$ is ``$\Delta p(n) = 1$", and $c$ is ``$n$ is a $q$". It follows that $\neg a$ is ``$p(n) + 1$ is a $q$", $\neg b$ is ``$\Delta p(n) = 2$", and $c$ is ``$n$ is a $p$". By our work for $\Delta p(n) \in \{1, 2\}$, we know $a \rightarrow b$ and $\neg a \rightarrow \neg b$. It suffices to show $a \rightarrow c$ and $\neg a \rightarrow \neg c$. In words, we want to show that ``If $p(n) + 1$ is a $p$, then $n$ is a $q$" and ``If $p(n) + 1$ is a $q$, then $n$ is a $p$". For convenience, consider our proof for $n \geq 2$ to avoid edge cases. Edge cases are left to the reader.

Suppose $p(n) + 1$ is a $p$. $p(n) - 1$ must be a $q$, say $q(l)$, as there cannot exist three consecutive $p$s. Similarly, we know $q(l + 1) = p(n) + 2$. In all, $p(n - 1) + 1 = q(l)$ ($q(l)$ is one more than the $p$ before $p(n)$), and $p(n + 1) + 1 = q(l + 1)$ (with similar reasoning). By the previous lemma, $q(l) = p(p(l)) + 1$ and $q(l + 1) = p(p(l + 1)) + 1$, so $p(p(l)) = p(n - 1)$ and $p(p(l + 1)) = p(n + 1)$. From $p$ and $q$ being strictly increasing, $n - 1 = p(l)$ and $n + 1 = p(l + 1)$. That is, $n - 1$ and $n + 1$ are $p$s, from which we must have that $n$ is a $q$.

Now suppose $p(n) + 1$ is a $q$, say $q(l)$. From the previous lemma again, $q(l) = p(p(l)) + 1$, so $p(p(l)) = p(n)$. Hence $p(l) = n$ as $p$ is increasing. In other words, $n$ is a $p$, as desired.}

If we know $\Delta p(n)$ for all $n \geq 1$, it is natural to think that we can write $p(n)$ recursively of itself. The following corollary confirms our intuition. \\

\corollary{$p(n + 1)$ is the number of $i$ in the range $1 \leq i \leq n$ appearing in $\{ p(1), p(2), p(3), \dots , p(n) \}$, added to $n + 1$.}{This is a straightforward induction argument following from the previous lemma and $p(n) \geq n$ (as $p$ is strictly increasing).}

A big idea here is that we have found a recurrence for $p(n)$ independent of $q(n)$, which we will exploit later. This statement is very powerful as the previous two-sequence definition of $p(n)$ and $q(n)$ was too clumsy to use effectively. Another other big idea is the methodology of using clever set construction to compute $p(p(n))$ (and therefore $q(p(n))=p(p(n))+p(n)$). We are very motivated to develop forms for $p(q(n))$ and $q(q(n))$, which we achieve in next few lemmas in the next section. A couple of corollaries follow from these new forms. The next section is not our main focus. Instead it shows the application of our crux method of attack and brings out the beauty of the coupled sequences $p$ and $q$. It is possible to skip the next section to see the relationship between $p$, $q$, and $\phi$.

\section{The Beauty of $p$ and $q$}
We will apply our method of set construction to discover beautiful relations in $p$ and $q$.\\

\corollary{$q(p(n)) = p(n) + q(n) - 1$.}{Follows from $p(p(n)) = q(n) + 1$ and the definition of $q$.}

\lemma{$p(q(n)) = p(n) + q(n)$.}{Consider  $[q(p(n)) + 1]$. There are $p(n)$ $q$s in $[q(p(n)) + 1]$, so there are $q(p(n)) + 1 - p(n) = p(p(n)) + 1$ $q$s. We know $q(n) = p(p(n)) + 1$, so there are $q(n)$ $p$s in $[q(p(n)) + 1]$. There cannot exist two consecutive $q$s, so $q(p(n)) + 1$ is a $p$, in fact the largest $p$ of $[q(p(n)) + 1]$. Hence $p(p(n)) = q(p(n)) + 1$. The previous corollary finishes.}

\corollary{If $\pair{p(n)}{q(n)}$ is a $p-q$ pair, then so is $\pair{p(n) + q(n)}{p(n) + 2 q(n)}$.}{We know $p(n) + q(n) = p(q(n))$ and $p(n) + 2q(n) = p(q(n)) + q(n) = q(q(n))$, as needed.}

The previous corollary reminds us greatly of the Fibonacci sequence. We finish with the most beautiful equation. \\

\corollary{$p(q(n)) = q(p(n)) + 1$.}{Follows from our forms for $p(q(n))$ and $q(p(n))$.}

\section{Linking $p$, $q$, and $\phi$}
Recall that we had the really motivating conjecture of $p(n) = \floor{\phi n}$ and $q(n) = \floor{\phi ^2 n}$. With our knowledge of the behavior of $p$ and $q$, we are finally ready to attack this conjecture. Unfortunately, we were unable to show equality. However, we present promising partial work. \\

\lemma{
Define by $E(n) = p(n) - \floor{n\phi}$. Then $E(n) \in \{ -1, 0, 1 \}$ for all $n \isint$.
}
{
We proceed by strong induction on $n$.\\

\induction{
By direct computation, we have $E(1) = E(2) = E(3) = E(4) = E(5) = E(6) = 0 \in \{ -1, 0, 1 \}$.
}
{
Assume that the result holds for $n = 1, 2, 3, \dots , k$ for some $k \geq 6$. We wish to show that the result holds for $n = k + 1$. That is, $E(k + 1) \in \{ -1, 0, 1 \}$. \\

There is a unique $j \isint$ for which $p(j) \leq k < p(j + 1)$. Our subgoal is to show that $j \leq k - 1$, so we can apply the inductive hypothesis to $p(j)$ and $p(j + 1)$. Denote by $j'$ the maximum positive integer for which $j' \leq k - 1 \text{   } (1)$ and $j' < \frac{k - E(j) + 1}{\phi} \text{   }(2)$. We can write $(2)$ equivalently as:
\[j'\phi < k - E(j') + 1 \Longleftrightarrow \floor{j'\phi} \leq k - E(j') \Longleftrightarrow \floor{j\phi} + E(j') \leq k \Longleftrightarrow p(j') \leq k.\]
Now $(1)$ is weaker than $(2)$ (as for $k \geq 6$, $\frac{k - E(j') + 1}{\phi} \leq \frac{k + 2}{\phi} \leq k - 1$), so we in fact have $j' = j$. We can now safely apply the induction hypothesis to $j$ and  $j+1$.

We know $j < \frac{k - E(j) + 1}{\phi} \leq \frac{k + 2}{\phi}$, so $j \leq \floor{\frac{k + 2}{\phi}}$. From $k < p(j + 1)$, we have
\begin{align*}
& k < \floor{(j + 1)\phi} + E(j + 1) \\
&\Longleftrightarrow k - E(j + 1) < \floor{(j + 1) \phi} \\
&\Longleftrightarrow k - E(j + 1) + 1 \leq \floor{(j + 1) \phi} \\
&\Longleftrightarrow k - E(j + 1) + 1 < (j + 1) \phi \\
&\Longleftrightarrow \frac{k - E(j + 1) + 1}{\phi} - 1 < j \\
&\Longrightarrow \frac{k - 1 + 1}{\phi} - 1 < j \\
&\Longrightarrow \floor{\frac{k}{\phi}} \leq j.
\end{align*}
This means that $\floor{\frac{k}{\phi}} \leq j \leq \floor{\frac{k + 2}{\phi}}$.\\

We now use the properties that for all $a, b \in \mathbb{R}$, we have $\floor{a + b} \leq \floor{a} + \floor{b} + 1$ and $\floor{a + b} \geq \floor{a} + \floor{b}$. Let $X = j - \floor{\frac{k}{\phi}} \geq 0$. Then 
\[X \leq \floor{\frac{k + 2}{\phi}} - \floor{\frac{k}{\phi}} \leq \floor{\frac{k}{\phi}} + \floor{\frac{2}{\phi}} + 1 - \floor{\frac{k}{\phi}} = 2.\]
Somewhat similarly, define $ Y = \floor{(k + 1) \phi} - \floor{k \phi}$ for some $Y \in \mathbb{Z}_{\geq 0}$. Again using the above floor property, we find
\[Y \leq \floor{k \phi} + \floor{\phi} + 1 - \floor{k \phi} = 2.\]
Using the other floor property, we find
\[Y \geq \floor{k \phi} + \floor{\phi} - \floor{k \phi} = 1.\]
In all, we know $X \in \{ 0, 1, 2 \}$ and $Y \in \{ 1, 2 \}$. \\

Note that the number of $l \leq k$ for which $l$ is in $p(1), p(2), p(3), \dots , p(n)$ is precisely $j$. It follows by Corollary 3, the recursive definition of $p(n)$ using itself, that $p(k + 1) = k + 1 + j$, and that:

\begin{align*}
E(k + 1) &= p(k + 1) - \floor{(k + 1) \phi} \\
&\stackrel{1}{=}k + 1 + \floor{k/ \phi} + X - \floor{k \phi} - Y \\
&\stackrel{2}{=} k + 1 + \floor{k(\phi - 1)} + X - \floor{k \phi} - Y \\
&= k + 1 + \floor{k \phi} - k + X - \floor{k \phi} - Y \\
&= 1 + X - Y.
\end{align*}
Note that step $1$ is from the definitions of $X$ and $Y$ and step $2$ is from a well known property of $\phi$. It remains to show that $1 + X - Y \in \{ -1, 0, 1 \}$. \\

The only case when $E(k + 1) \not \in \{ -1, 0, 1 \}$ is when $(X, Y) = (2, 1)$. It remains for us to show that this case is impossible. Suppose for the sake of contradiction that $E(j) \leq 0$. Then from $j < \frac{k - E(j) + 1}{\phi}$, we have $j < \frac{k + 1}{\phi}$, so $j \leq \floor{\frac{k}{\phi}} + 1$ and hence $X \leq 1$, a contradiction. This means $E(j) = -1$ and $j = \frac{k + 2}{\phi}$, from which we find $\floor{\frac{k + 2}{\phi}} = \floor{\frac{k}{\phi}} + 2$ as $X = 2$. We rewrite this equation as follows:

\begin{align*}
&\floor{\frac{k + 2}{\phi}} = \floor{\frac{k}{\phi}} + 2 \\
&\Longleftrightarrow \floor{(k + 2)(\phi - 1)} = \floor{k (\phi - 1)} + 2 \\
&\Longleftrightarrow \floor{(k + 2 )\phi} - k -2 = \floor{k \phi} - k + 2 \\
&\Longleftrightarrow \floor{(k + 2) \phi} = \floor{k \phi} + 4 \text{   } (3).
\end{align*}

Using the fact that $Y = 1$, or that $\floor{(k + 1) \phi} = \floor{k\phi} + 1$ with $(3)$, we get $\floor{(k + 2) \phi} = \floor{(k + 1)\phi} + 3$, which is impossible. \\

We have therefore shown $E(k + 1) \in \{ -1, 0, 1 \}$ and hence completed the strong induction.
}
}
\indent The main takeaway from this proof is that we used local information relative to the functions $f(n) = \floor{n \phi}$  and $g(n) = \floor{n / \phi}$ by applying it to what we already know about $p(n)$. If we can extract stronger local information of $f$ and $g$, we can possibly sharpen our result to $p(n) = \floor{n \phi}$. Regardless of potential improvement, our result is already very strong. If we are given any $n$, we can guarantee $p(n)$ to be within $1$ away from $\floor{ n\phi}$ without having to perform any recursion.

\section{Future Work}
Our method of set construction has weak application to the primes. Let $P$ denote the sequence of all primes in increasing order. That is, $P$ contains $2, 3, 5, 7, \dots $ in that order. Define by $Q$ the gaps of composites that $P$ left. That is, $Q$ contains $4, 6, 8, 9, 10, 12, 14, 15, 16, 18, 20, 21, 22, \dots$ in that order. We claim that $Q(P(n) - n - 1) = P(n) - 1$ for $n \geq 3$. Consider the set $[P(n)]$. There are $n$ $P$s in this set, so there are $P(n) - n - 1$ $Q$s in this set (as we do not count $1$ in either $P$ or $Q$). For $n \geq 3$, $P(n) - 1$ must be a $Q$ as the only consecutive primes are $2$ and $3$. Hence $Q(P(n) - n - 1) = P(n) - 1$, as desired. We believe that it is possible to create constructions that better exploit the properties of the primes and create more useful equations. Here we only used the fact that most primes are nonconsecutive.

\section{Conclusion}
We used elementary techniques to prove non-elementary results. Our initial confusing definition of $p$ and $q$ gave the unexpected conjecture of $p(n) = \floor{\phi n}$ and $q(n) = \floor{\phi^2 n}$. We were close to proving this relationship by showing that $p(n)$ is at most $1$ away from $\floor{\phi n}$. Our elementary technique of clever set constructions allowed us to convert qualitative data about $p$ and $q$ into quantitative data.

\printbibliography

\end{document}